\documentclass[12pt,reqno]{amsart}
\usepackage[cp1251]{inputenc}
\usepackage{amsfonts}
\usepackage{amssymb,amsmath,amsthm,eufrak,euscript}
\usepackage{graphicx,fancybox}
\newtheorem{lemma}{Lemma}
\newtheorem{theorem}{Theorem}

\newtheorem{example}{Example}
\newtheorem{remark}{Remark}
\numberwithin{equation}{section}

\begin{document}

\sloppy

\begin{center}
{\bf \Large Rational integrals of 2-dimensional \\ geodesic flows: new examples}
\end{center}

\medskip

\begin{center}
{\bf Sergei Agapov , Vladislav Shubin}
\footnote{Both authors are supported by the Mathematical Center in Akademgorodok under the agreement No. 075-15-2019-1675 with the Ministry of Science and Higher Education of the Russian Federation.}
\end{center}

\medskip

\begin{quote}
\noindent{\bf Abstract. }{\small This paper is devoted to searching for Riemannian metrics on 2-surfaces whose geodesic flows admit a rational in momenta first integral with a linear numerator and denominator. The explicit examples of metrics and such integrals are constructed. Few superintegrable systems are found having both a polynomial and a rational integrals which are functionally independent of the Hamiltonian.
}

\medskip


\end{quote}

\medskip

\section{Introduction}

A geodesic flow of the Riemannian metric $ds^2=g_{ij}dx^idx^j$ on a 2-surface is given by the Hamiltonian system
$$
\dot{x}^j = \{x^j,H\}, \quad \dot{p}_j = \{p_j,H\}, \quad H = \frac{1}{2} g^{ij} p_ip_j, \quad i,j=1,2, \eqno(1.1)
$$
the Poisson bracket has the form:
$$
\{F,H\} = \sum_{i=1}^2 \left ( \frac{\partial F}{\partial x^i} \frac{\partial H}{\partial p_i} - \frac{\partial F}{\partial p_i} \frac{\partial H}{\partial x^i} \right ).
$$
The geodesic flow (1.1) is called completely integrable, if there exists an additional first integral, i.e. a function $F(x,p)$ such that $\frac{dF}{dt} = \{F,H\} \equiv 0,$ and $F, H$ are functionally independent a.e..

Searching for metrics on 2-surfaces with integrable geodesic flows is a classical problem. There exist certain topological obstacles to the global integrability: on surfaces of any genus larger than~1 there are no analytic Riemannian metrics with analytically integrable geodesic flows (see~\cite{1}). Note that in the smooth case this result may be wrong in general (see~\cite{2}, Chapter 3).

In the overwhelming majority of known examples the first integrals are polynomials in momenta. The polynomial integrals of small degrees are well-studied (see, e.g.,~\cite{3},~\cite{4}): the existence of a linear integral is related to the presence of a cyclic variable; the quadratic one is related to a possibility of the separation of variables. It was proved in~\cite{2},~\cite{5} that on 2-surfaces there exist Riemannian metrics whose geodesic flow admits a local polynomial integral of an arbitrarily high degree which is independent of the Hamiltonian. Rational in momenta first integrals of mechanical systems are also of an interest, they have been studied in many papers (see, e.g.,~\cite{6}---~\cite{14}). Apparently, J.G. Darboux was one of the first who started an investigation of such integrals for geodesic flows (\cite{15}). It was proved in~\cite{16} that on 2-surfaces there exist analytic Riemannian metrics with an additional local rational integral with a numerator and a denominator of arbitrarily high degrees. The proof is based on applying the Cauchy-Kovalevskaya theorem. However, the explicit construction of such metrics and rational integrals (as well as polynomial integrals of high degrees) faces with various difficulties: one can find few known examples in~\cite{17}---~\cite{19} (see also~\cite{14}). In this paper we investigate this problem (partially following the ideas outlined in~\cite{15}) in the simplest case when the degrees of a numerator and a denominator of the rational integral are equal to 1.

\section{Reduction to a sole equation}

Suppose that the geodesic flow of the metric $ds^2 = \Lambda(x,y) (dx^2+dy^2)$ on a 2-surface with coordinates $x,y$ admits a rational integral
$$
F=\frac{f_1(x,y)p_1+g_1(x,y)p_2}{f_2(x,y)p_1+g_2(x,y)p_2}.
$$
The Hamiltonian has the form $H = \frac{p_1^2+p_2^2}{2 \Lambda}.$ The Poisson bracket vanishes which produces
\begin{gather}
2 \Lambda (f_2 {f_1}_x - f_1 {f_2}_x) + \Lambda_y (f_2g_1-f_1g_2) = 0, \\
2 \Lambda (g_2 {f_1}_x - g_1 {f_2}_x + f_2 ({f_1}_y + {g_1}_x) - f_1 ({f_2}_y+{g_2}_x)) + \Lambda_x (f_1g_2-f_2g_1) = 0, \\
2 \Lambda (f_2 {g_1}_y - f_1 {g_2}_y + g_2 ({f_1}_y + {g_1}_x) - g_1 ({f_2}_y+{g_2}_x)) - \Lambda_y (f_1g_2-f_2g_1) = 0, \\
2 \Lambda (g_2 {g_1}_y - g_1 {g_2}_y) + \Lambda_x (f_1g_2-f_2g_1) = 0.
\end{gather}
In slightly different notations this system was written down in~\cite{16}, an existence of analytic solutions (generally speaking, local ones) was also established there. In this paper we construct such solutions in an explicit form. According to~\cite{14} (2.1) --- (2.4) imply
\begin{equation}
\left( \frac{f_1 - i g_1}{f_2 - i g_2} \right)_x - i \left( \frac{f_1 - i g_1}{f_2 - i g_2} \right)_y = 0.
\end{equation}
Let us rewrite this equality in the following way
$$
\left( u - i v \right)_x - i \left( u - i v \right)_y = 0,
$$
where
$$
u = \frac{f_1 f_2 + g_1 g_2}{f_2^2+g_2^2}, \quad v = \frac{f_2 g_1 - f_1 g_2}{f_2^2+g_2^2}.
$$
Consequently, $u_{xx} + u_{yy} = v_{xx} + v_{yy} = 0,$ i.e. $u(x,y), v(x,y)$ are two conjugate harmonic functions:
\begin{equation}
u_x = v_y, \qquad u_y = - v_x.
\end{equation}
Notice that if $v(x,y) \equiv 0,$ then the first integral $F$ does not depend on momenta.
We may find $f_1(x,y), g_1(x,y)$ in the following way:
\begin{equation}
f_1 = f_2 u - g_2 v, \quad g_1 = g_2 u + f_2 v.
\end{equation}
Denote $\Lambda(x,y) = e^{\lambda(x,y)}.$
Due to (2.6), (2.7) the system (2.1) --- (2.4) is equivalent to the following two equations:
\begin{gather}
(f_2^2+g_2^2)v\lambda_y -2 f_2 (v {g_2})_x + 2 f_2^2 v_y + 2 g_2 v {f_2}_x = 0, \\
(f_2^2+g_2^2)v\lambda_x -2 g_2 (v {f_2})_y + 2 g_2^2 v_x + 2 f_2 v {g_2}_y = 0.
\end{gather}
We will restrict ourselves to searching for solutions such that $f_2(x,y)$ and $g_2(x,y)$ do not vanish simultaneously. Without loss of generality assume that $f_2^2 + g_2^2 = 1$. Put $f_2(x,y) = \sin\mu(x,y)$, $g_2(x,y) = \cos\mu(x,y)$, then (2.8), (2.9) imply that
\begin{gather}
v (\lambda_y + 2 \mu_x) - 2 \sin\mu\,(v_x \cos\mu - v_y \sin\mu) = 0, \\
v (2 \mu_y - \lambda_x) - 2 \cos\mu\,(v_x \cos\mu - v_y \sin\mu) = 0.
\end{gather}
Multiply (2.10) by $\sin\mu$, (2.11) by $\cos\mu$ and sum them up; the obtained equality is equivalent to:
\begin{equation}
\left(v e^{\frac{\lambda}{2}}\cos\mu\right)_x - \left(v e^{\frac{\lambda}{2}}\sin\mu\right)_y = 0.
\end{equation}
By the analogous way, multiplying (2.10) by $\cos\mu$, (2.11) by $\sin\mu$ and subtracting one from the other, one obtains
\begin{equation}
\left(e^{\frac{\lambda}{2}}\cos\mu\right)_y + \left(e^{\frac{\lambda}{2}}\sin\mu\right)_x = 0.
\end{equation}
Put $\sigma = e^{\frac{\lambda}{2}}\cos\mu$, $\tau = -e^{\frac{\lambda}{2}}\sin\mu$. The equations (2.12), (2.13) take the form
\begin{equation}
\notag
(v\sigma)_x + (v\tau)_y = 0,\quad \sigma_y = \tau_x.
\end{equation}
The equality $\sigma_y = \tau_x$ implies that there exists a function $\psi(x,y)$ such that $\psi_x = \sigma$, $\psi_y = \tau$. Thus, the first equation is equivalent to
\begin{equation}
(v\psi_x)_x + (v\psi_y)_y = 0.
\end{equation}
Let us express $f_2$, $g_2$ in terms of $\psi$. In view of the substitutions made above we have $f_2 = -e^{-\frac{\lambda}{2}}\psi_y, \quad g_2 = e^{-\frac{\lambda}{2}}\psi_x.$ Since the coefficients of $F$ are defined up to a multiplication by an arbitrary function, due to (2.7) one may assume that $f_2 = -\psi_y$, $g_2 = \psi_x.$ Expressing the other unknown functions in terms of $\psi,$ we obtain
\begin{gather}
f_2 = -\psi_y, \qquad g_2 = \psi_x, \qquad f_1 = -\psi_y u - \psi_x v, \qquad g_1 = \psi_x u - \psi_y v,\\
\lambda = \ln(\psi_x^2 + \psi_y^2), \qquad \Lambda = \psi_x^2 + \psi_y^2.
\end{gather}

Thus the following theorem holds true.
\begin{theorem}
Let $u(x,y), v(x,y)$ be two arbitrary conjugate harmonic functions satisfying (2.6). Also assume that $\psi(x,y)$ is an arbitrary function satisfying (2.14) and functions $f_1(x,y)$, $g_1(x,y)$, $f_2(x,y)$, $g_2(x,y)$ and $\Lambda(x,y)$ satisfy (2.15) and (2.16).
Then the geodesic flow of the metric $ds^2 = \Lambda(x,y)(dx^2 + dy^2)$ admits a rational integral
\begin{equation}
\notag
F = \frac{f_1(x,y)p_1 + g_1(x,y)p_2}{f_2(x,y)p_1 + g_2(x,y)p_2}.
\end{equation}
\end{theorem}

Summing it up, given a solution $\psi(x,y)$ to the equation (2.14), one can find the metric and all the coefficients of the first integral $F$ using (2.15), (2.16).

\section{Construction of solutions to the key equation}

Let us make the change of variables $(x,y) \rightarrow (u(x,y), v(x,y))$ in (2.14). If the gradient of $v(x,y)$ does not vanish, then this change is non-degenerate and (2.14) is equivalent to
\begin{equation}
\varphi_{uu} + \varphi_{vv} + \frac{1}{v}\varphi_v = 0,
\end{equation}
where $\psi(x,y) = \varphi\big(u(x,y), v(x,y)\big)$.
In case of a degenerate change the question about an equivalence of these equations requires an additional investigation, however, by the direct differentiation it is easy to establish that the following statement holds true.
\begin{lemma}
Let $\varphi(u, v)$ be a solution to (3.1). Then for any pair of conjugate harmonic functions $u(x,y), v(x,y)$ satisfying (2.6), the function
\begin{equation}
\notag
\psi(x,y) = \varphi\big(u(x,y), v(x,y)\big)
\end{equation}
is a solution to (2.14).
\end{lemma}

We shall search for complex-valued solutions to (3.1) via the method of separations of variables. Substituting $\varphi(u,v)=\alpha(u)\beta(v)$ into (3.1) and dividing by $\alpha(u)\beta(v)$, one obtains
\begin{equation}
\notag
\frac{\alpha''(u)}{\alpha(u)} = -\frac{\beta''(v)}{\beta(v)}-\frac{\beta'(v)}{v\beta(v)}.
\end{equation}
Since the left-hand side of this equality does not depend on $v$, the right-hand side does not depend on $u$, then both sides are constants. So we obtain
\begin{equation}
\notag
\frac{\alpha''(u)}{\alpha(u)} = \nu,\quad \frac{\beta''(v)}{\beta(v)}+\frac{\beta'(v)}{v\beta(v)}=-\nu
\end{equation}
for a certain complex constant $\nu$ or, equivalently,
\begin{gather}
\alpha''(u) = \nu \alpha(u),\\
\beta''(v) + \frac{1}{v}\beta'(v) + \nu\beta(v)=0.
\end{gather}
Thus, if the functions $\alpha(u)$, $\beta(v)$ satisfy (3.2), (3.3), then the function $\varphi(u,v)=\alpha(u)\beta(v)$ is the solution to (3.1).

For convenience put $s = \sqrt{\nu}$. The general solution to (3.2) has the form $\alpha(u) = a_1 e^{su} + a_2 e^{-su}$, where $a_1$, $a_2$ are arbitrary complex numbers. The equation (3.3) after multiplying by $v^2$ and making the change of variables $\tilde{v} = s v$ is equivalent to
\begin{equation}
\notag
\tilde{v}^2 \tilde{\beta}'' + \tilde{v}\tilde{\beta}' + \tilde{v}^2\tilde{\beta}=0,
\end{equation}
where $\tilde{\beta}(\tilde{v}) = \tilde{\beta}(sv) = \beta(v).$ The obtained equation is the Bessel equation (\cite{20}) and has the general solution of the form
\begin{equation}
\notag
\tilde{\beta}(\tilde{v}) = b_1 J_0(\tilde{v}) + b_2 Y_0(\tilde{v}),
\end{equation}
where $b_1$, $b_2$ are arbitrary complex numbers, $J_0$ is the Bessel function of the first kind, $Y_0$ is the Bessel function of the second kind (\cite{20}), i.e.
\begin{equation}
\notag
\beta(v) = b_1 J_0(sv) + b_2 Y_0(sv).
\end{equation}
It will be convenient for us to write down the general solution in a slightly different form. It is known (\cite{20}) that the function $Y_0(z)$ can be represented in the form
\begin{equation}
Y_0(z) = \frac{2}{\pi}\ln\Big(\frac{z}{2}\Big)J_0(z) + \sum_{k=0}^\infty d_k z^{2k},
\end{equation}
where coefficients $d_k$ are certain real numbers (one can find their explicit form in~\cite{20}), and the series $\sum_{k=0}^\infty d_k z^{2k}$ converges on the whole plane of the complex variable $z$. Notice that
\begin{equation}
\notag
\beta_1(s,v) = Y_0(sv)-\frac{2}{\pi}\ln(s)J_0(sv)
\end{equation}
is the analytic function of the complex variable $s$ on the whole plane for any real $v \neq 0$ and simultaneously is the partial solution to (3.3) which is linearly independent of $J_0(sv)$ for any complex $s.$ Consequently, the function
\begin{equation}
\notag
\beta(v) = b_1 J_0(sv) + b_2 \beta_1(s, v)
\end{equation}
with arbitrary complex constants $b_1$, $b_2$ is the general solution to (3.3).

Thus, the function
\begin{equation}
\varphi(u,v,s) = \alpha(u)\beta(v) = c_1 e^{su}J_0(sv) + c_2 e^{-su}J_0(sv) + c_3 e^{su} \beta_1(s, v) + c_4 e^{-su} \beta_1(s, v),
\end{equation}
is the solution to the equation (3.1) for  any $s\in\mathbb{C},$ where $c_1, c_2, c_3, c_4$ are arbitrary complex numbers. Due to the linearity of (3.1) we obtain that any linear combination of the functions of the form (3.5) with different $s$ is also a solution to (3.1). Moreover, notice that due to the parity of $J_0$ and due to (3.4) we have
\begin{equation}
\notag
\beta_1(-s, y) = \beta_1(s, y) + i\pi J_0(sv).
\end{equation}
Consequently, the function of the form (3.5) is a linear combination of the functions of the form
\begin{equation}
\varphi(u,v,s) = c_1 e^{su} J_0(sv) + c_2 e^{su} \beta_1(s, v),
\end{equation}
with opposite values of $s$. Here $c_1$, $c_2$ are arbitrary complex numbers.

Notice that since the classical integral is the limit of certain Riemannian sums consisting of finite number of summands, then one may assume that an integral of functions of the form (3.6) over $s$ is also a solution. Let us prove the following statement.

\begin{lemma}
Assume that $n$, $m$ are arbitrary natural numbers, $\gamma_{1,1}$,...,$\gamma_{1,n}$, $\gamma_{2,1}$,...,$\gamma_{2,m}$ are arbitrary rectifiable (i.e. having a finite length) curves (closed or unclosed) on the complex plane. Also let $a_1(s)$,...,$a_n(s)$, $b_1(s)$,...,$b_m(s)$ be bounded functions defined on the curves $\gamma_{1,1}$,...,$\gamma_{1,n}$, $\gamma_{2,1}$,...,$\gamma_{2,m}$ accordingly. Then the function
\begin{equation}
\varphi(u, v) = \sum_{j=1}^n \int_{\gamma_{1,j}} a_j(s) e^{su}J_0(sv)ds +
\sum_{j=1}^m\int_{\gamma_{2,j}} b_j(s) e^{su}\Big(Y_0(sv) - \frac{2}{\pi}\ln(s)J_0(sv)\Big)ds
\end{equation}
is the solution to the equation (3.1), which is analytic in both real variables $u$, $v$ at any point where $v\neq 0$. If $b_j(s) \equiv 0, \ 1 \leqslant j \leqslant m,$ then this solution is analytic on the whole plane of the real variables $u, v.$
\end{lemma}

\begin{proof}
Due to the linearity of the equation (3.1) it is enough to prove that each of the summands of (3.7) is the solution to this equation. We have
\begin{gather}
\notag
\frac{\partial^2}{\partial u^2} \big(a_1(s)e^{su}J_0(sv)\big) = s^2 a_1(s)e^{su}J_0(sv),\\
\notag
\frac{\partial}{\partial v} \big(a_1(s)e^{su}J_0(sv)\big) = s a_1(s)e^{su}J_0'(sv),\\
\notag
\frac{\partial^2}{\partial v^2} \big(a_1(s)e^{su}J_0(sv)\big) = s^2 a_1(s)e^{su}J_0''(sv).
\end{gather}

The function $a_1(s)$ and the curve $\gamma_{1,1}$ are bounded, so due to the Lebesgue theorem the derivative of the integral coincides with the integral of the derivative, i.e.
\begin{gather}
\notag
\frac{\partial^2}{\partial u^2} \int_{\gamma_{1,1}}a_1(s)e^{su}J_0(sv)ds = \int_{\gamma_{1,1}} s^2 a_1(s)e^{su}J_0(sv)ds,\\
\notag
\frac{\partial}{\partial v} \int_{\gamma_{1,1}} a_1(s)e^{su}J_0(sv)ds = \int_{\gamma_{1,1}} s a_1(s)e^{su}J_0'(sv)ds,\\
\notag
\frac{\partial^2}{\partial v^2} \int_{\gamma_{1,1}} a_1(s)e^{su}J_0(sv)ds = \int_{\gamma_{1,1}} s^2 a_1(s)e^{su}J_0''(sv)ds.
\end{gather}
It is left to substitute the obtained equalities into (3.1) and to use the fact that $J_0$ satisfies the Bessel equation. In an analogous way it can be proved that the other summands also satisfy the equation (3.1).

Let us prove now that the obtained solution is analytic at $v\neq 0$. For example, let us prove the analyticity of the first summand of the second sum in (3.7) at the point $(u_0, v_0)$. Denote
\begin{equation}
\notag
A(u, v) = \int_{\gamma_{2,1}} b_1(s)e^{su}\left(Y_0(sv)-\frac{2}{\pi}\ln(s)J_0(sv)\right)ds.
\end{equation}
Let $|u-u_0|<\varepsilon$, $|v-v_0|<\varepsilon$ for some small enough $\varepsilon > 0$ (in particular $\varepsilon < |v_0|$), and let $\Gamma$ be an arbitrary contour laying in the complex neighborhood $|v-v_0|<\varepsilon$. We have
\begin{equation}
\notag
\int_\Gamma A(u, v)dv =
\int_{\gamma_{2,1}} b_1(s)e^{su} ds \int_\Gamma \left(Y_0(sv)-\frac{2}{\pi}\ln(s)J_0(sv)\right) dv.
\end{equation}
Changing the order of integration is possible due to the fact that the integrand is a continuous and bounded function on the set $\gamma_{2,1}\times\{u\mid|u-u_0|<\varepsilon\}\times\Gamma$. Due to (3.4) we have
\begin{equation}
\notag
\beta_1(s, v) = Y_0(sv)-\frac{2}{\pi}\ln(s)J_0(sv) = \frac{2}{\pi}\ln\left(\frac{v}{2}\right)J_0(sv)+\sum_{k=0}^\infty d_k s^{2k}v^{2k}.
\end{equation}
Since the variable $v$ does not vanish in the neighborhood $|v-v_0|<\varepsilon$ for small enough $\varepsilon$, then $\beta_1(s, v)$ is an analytic function in $v$ in this neighborhood. Then
\begin{equation}
\notag
\int_\Gamma \left(Y_0(sv)-\frac{2}{\pi}\ln(s)J_0(sv)\right) dv = 0
\end{equation}
holds true due to the Cauchy theorem. Thus, the Morera's theorem (\cite{21}) implies that $A(u, v)$ is the analytic function in $v$ on $|v-v_0|<\varepsilon$. The analyticity of this function in $u$ on $|u-u_0|<\varepsilon$ can be proved analogously. Further, the Hartog's theorem (\cite{21}) implies that $A(u, v)$ is an analytical function in both variables $(u, v)$ on the set $|u-u_0|<\varepsilon$, $|v-v_0|<\varepsilon$. Since $(u_0, v_0)$ is an arbitrary point such that $v_0\neq 0$, then the first summand of the second sum is analytic in both variables $u$ and $v$ everywhere except the points where $v = 0$. The analyticity of the other summands can be proved analogously. Moreover, the arguments mentioned above hold true for the summands of the first sum for any point $(u_0, v_0)$.
\end{proof}

\begin{remark}
Generally speaking, the formula (3.7) gives complex-valued solutions. But since the equation (3.1) is linear and real-valued, then the real and imaginary parts of (3.7) are also solutions. Moreover, one can obtain the imaginary part of (3.7) considering the real part of a function of the same form, where functions $a_j(s)$, $b_j(s)$ are multiplied by $-i.$ Since $a_j(s), b_j(s)$ are arbitrary functions, then it is enough to consider only a real part of the functions of the form (3.7) to obtain real-valued solutions to the equation (3.1).
\end{remark}

\begin{remark}
The solutions of the form (3.7) include the ones which are linear combinations of the solutions of the form (3.6). To demonstrate this, it is enough to put
\begin{equation}
\notag
a_1(s) = \frac{1}{2\pi i}\sum_{j=1}^n \frac{c_{1,j}}{s-s_{1,j}},\quad
b_1(s) = \frac{1}{2\pi i}\sum_{j=1}^m \frac{c_{2,j}}{s-s_{2,j}},
\end{equation}
where $s_{1,j}$, $s_{2,j}$, $c_{1,j}$, $c_{2,j}$ are some complex numbers, and $\gamma_{1,1},$ $\gamma_{2,1}$ are closed contours surrounding all the values $s_{1,j},$ $s_{2,j}$ accordingly. Substituting these functions and contours in (3.7) we obtain the solution
\begin{equation}
\varphi(u, v) = \sum_{j=1}^n c_{1,j} e^{s_{1,j} u}J_0(s_{1,j} v)+
\sum_{j=1}^m c_{2,j} e^{s_{2,j} u}\Big(Y_0(s_{2,j} v)-\frac{\pi}{2}\ln(s_{2,j})J_0(s_{2,j} v)\Big).
\end{equation}
Moreover, put $a_1(s) = \frac{m!}{2\pi i}\frac{1}{(s-s_0)^{m+1}}$ and assume that the contour $\gamma_{1,1}$ is a circle of an arbitrary radius with the center in the point $s_0$, then due to the formula of derivatives of the Cauchy integral we obtain the following solution
\begin{equation}
\varphi(u, v) = \frac{d^m}{ds^m}\big(e^{su}J_0(sv)\big)\Big|_{s=s_0}.
\end{equation}

Thus, the formula (3.7) includes some partial cases which are linear combinations (for different $s$) of  derivatives of the solutions of the form (3.6) over $s$ of an arbitrary order $m \in \mathbb{N}.$
\end{remark}

\section{Examples and the main result}

Remark 2 allows one to construct various examples of metrics and rational integrals. Let us show a couple of examples.

\begin{example}
Assume that $u = x$, $v = y$. The function $\varphi(u, v) = e^u J_0(v)$ has the form (3.8) at $n=1$, $m=0$, $c_{1,1} = 1$, $s_{1,1} = 1$. Using the equality $J_0'(y) = -J_1(y)$ (see~\cite{20}) we obtain the metric
\begin{equation}
\notag
\Lambda(x,y) = e^{2x}(J_0^2(y)+J_1^2(y)).
\end{equation}
This is the analytic function on the whole plane of real variables $(x,y)$ and positive everywhere (due to the properties of zeros of the Bessel function (see~\cite{20})). The geodesic flow of this metric admits the first integral
\begin{equation}
\notag
F = x - y\frac{J_0(y)p_1 - J_1(y)p_2}{J_1(y)p_1 + J_0(y)p_2}.
\end{equation}
\end{example}

\begin{example}
Consider the solution (3.9) at $m=2$, $s_0 = 0$, i.e. $\varphi(u, v) = u^2-\frac{1}{2}v^2$. Put $u=e^x\cos(y)$, $v=e^x\sin(y)$. Then
\begin{equation}
\Lambda(x,y) = e^{4x}(1 + 3\cos^2(y)).
\end{equation}
As in the previous example, this is the analytic function on the whole plane of real variables $(x,y)$ and positive everywhere. In this case the geodesic flow admits the first integral
\begin{equation}
\notag
F = e^x\frac{\sin(y)p_1+2\cos(y)p_2}{3\sin(y)\cos(y)p_1+(2-3\sin^2(y))p_2}.
\end{equation}
Moreover, notice that the change $u=e^x\cos(y)$, $v=e^x\sin(y)$ separates the variables and produces the Liouville metric
\begin{equation}
\notag
ds^2 = (4 u^2 + v^2)(du^2+dv^2).
\end{equation}
Consequently, the geodesic flow of the metric (4.1) admits a quadratic integral $F_1,$ in canonical coordinates it has the form
\begin{equation}
\notag
F_1 = \frac{e^{-2 x}}{1+3\cos^2(y)}\Big(3\cos^2(y)\sin^2(y)p_1^2+\left(4 \cos^4(y)-\sin^4(y)\right)p_2^2\Big)+e^{-2 x}\sin(2y)p_1p_2.
\end{equation}
It is not difficult to check that the functions $H$, $F$, $F_1$ are functionally independent a.e. Thus we constructed the Riemannian metric which is analytic on the whole real plane $(x,y)$ such that its geodesic flow is superintegrable.
\end{example}

\begin{example}
Consider the solution (3.9) at $m=3$, $s_0 = 0$, for convenience let us multiply it by $\frac{2}{3},$ i.e. $\varphi(u, v) = \frac{2}{3}u^3-u\,v^2$. Also put $u=e^x\cos(y)$, $v=e^x\sin(y)$. We obtain the metric
\begin{equation}
\notag
\Lambda(x,y) = e^{6 x} \big(4\cos^4(y) + \sin^4(y)\big),
\end{equation}
which is also an analytic function on the whole plane $(x,y)$ and positive everywhere. The geodesic flow admits the first integral
\begin{equation}
\notag
F = 2 e^x \frac{2\sin(2 y)p_1 + (1+3\cos(2 y))p_2}{(\sin(y)+5 \sin(3 y))p_1 + (3 \cos(y)+5 \cos(3 y))p_2}.
\end{equation}
In this example the change $u=e^x\cos(y)$, $v=e^x\sin(y)$ also separates the variables and produces the Liouville metric
\begin{equation}
\notag
ds^2 = (4 u^4 + v^4)(du^2+dv^2),
\end{equation}
its geodesic flow admits the quadratic integral
\begin{equation}
\notag
F_1 = e^{-2 x} \frac{4 \cos^4(y) \big(p_2 \cos(y)+p_1 \sin(y)\big)^2-\sin^4(y) \big(p_1 \cos(y)-p_2 \sin(y)\big)^2}{4 \cos^4(y)+\sin^4(y)}.
\end{equation}
This is also the example of an analytic Riemannian metric on the whole real plane with the superintegrable geodesic flow.
\end{example}

\begin{example}
Construct a local metric which is expressed in terms of elementary functions of another kind.
Due to (3.4) we have $\beta_1(0, v) = \frac{2}{\pi}\ln(v) + d_0$ and $(\beta_1)_s(0,v) = 0$. Consequently the function
\begin{equation}
\notag
\tilde{\varphi}(u, v) = \frac{\pi}{2}\frac{d}{ds}\big(e^{su}\beta_1(s,v)\big)\Big|_{s=0} = u\ln(v) + \frac{\pi}{2}d_0u,
\end{equation}
is also the solution obtained via the same ideas as the solution of  the form (3.9). Since $\frac{\pi}{2}d_0u$ is also a solution, one can construct the metric and the first integral from the solution $\varphi(u, v) = u\ln(v)$. Put $u = x$, $v = y$. We obtain the metric
\begin{equation}
\notag
\Lambda(x, y) = \frac{x^2}{y^2} + \ln^2(y),
\end{equation}
which is the positive analytic function at $y > 1$. The geodesic flow admits the first integral
\begin{equation}
\notag
F = \frac{(x^2 + y^2\ln(y))p_1 - xy(\ln(y) - 1)p_2}{x\,p_1-y\ln(y)p_2}.
\end{equation}
\end{example}

Summing up the arguments mentioned above, let us formulate and prove the following

\begin{theorem}
Let $u(x,y)$, $v(x,y)$ be arbitrary conjugate harmonic functions satisfying the relation (2.6). Assume that $a_j(s)$, $b_l(s)$, $\gamma_{1,j},$ $\gamma_{2,l},$ where $1\leqslant j\leqslant n$, $1\leqslant l\leqslant m$ are functions and curves accordingly satisfying the conditions of Lemma 2. Also assume that $\varphi(u,v)$ has the form (3.7) and $\psi(x, y)=\mathrm{Re}\,\varphi\big(u(x, y), v(x, y)\big)$. Then the geodesic flow of the metric
\begin{equation}
ds^2 = \Lambda(x,y)(dx^2 + dy^2) , \qquad \Lambda(x, y) = \psi_x^2(x, y) + \psi_y^2(x, y)
\end{equation}
admits a rational first integral
\begin{equation}
F =
\frac{-\big(\psi_y u + \psi_x v\big)p_1 + \big(\psi_x u-\psi_y v\big)p_2}{-\psi_y p_1 + \psi_x p_2},
\end{equation}
wherein $\Lambda(x,y)$ is an analytic function at the points where $v(x,y)\neq 0$. Moreover, if $b_l(s) \equiv 0, \ l = 1, \ldots, m,$ then $\Lambda(x,y)$ is an analytic function on the whole plane $(x,y).$
\end{theorem}

\begin{proof}
Lemma 1, Lemma 2, and Remark 1 imply that $\psi(x,y)$ is a solution to the equation (2.14). Thus, the fact that (4.3) is the first integral of the geodesic flow of the metric (4.2) follows directly from the statement of Theorem 1. It is left to prove the analyticity of the function $\Lambda(x,y).$ This fact holds true since the harmonic functions are analytic ones in both real variables (\cite{21}) and the superpositions of analytic functions is also the analytic one. Indeed, if $v(x, y)\neq 0$ in a certain point, then in a certain neighborhood of this point $\psi(x, y)$ is an analytic function as the real part (for real $x$ and $y$) of the superposition of analytic functions $\varphi(u, v)$ and $\big(u(x, y), v(x, y)\big)$. Consequently, its derivatives are also analytic functions. If $b_l(s) \equiv 0$ for all $l$ and $s,$ then these arguments hold true for any point of the plane $(x, y)$.
\end{proof}

Let us show an example which was obtained directly via Theorem 2.
\begin{example}
Let $\gamma_{1,1}$ be a segment $[0,1]$ and assume that $a_1(s) = 1$, $u = x$, $v = y$. We obtain the function
\begin{equation}
\notag
\psi(x,y) = \int_0^1 e^{sx}J_0(sy)ds
\end{equation}
which produces the metric (due to (4.2))
\begin{equation}
\notag
\Lambda(x,y) = \Big(\int_0^1 s e^{sx} J_0(sy)ds\Big)^2+\Big(\int_0^1 s e^{sx} J_1(sy)ds\Big)^2.
\end{equation}
Since $\Lambda(0, 0) = 1/4 > 0,$ then due to the continuity $\Lambda(x,y)>0$ in a certain neighborhood of the origin. The geodesic flow of the given metric admits the first integral
\begin{equation}
\notag
F = x + y\frac{-p_1\int_0^1 s e^{sx}J_0(sy)ds + p_2\int_0^1 s e^{sx} J_1(sy)ds}{p_1\int_0^1 s e^{sx}J_1(sy)ds + p_2\int_0^1 s e^{sx} J_0(sy)ds}.
\end{equation}
\end{example}

\section{Conclusion}

Theorem 2 provides a general algorithm to construct any number of explicit examples of metrics whose geodesic flows admit rational in momenta first integral with a linear numerator and denominator. In general such metrics are defined and do not vanish only locally. However, as Example 1 shows, sometimes such metrics happen to be defined on the whole plane. Examples 2 and 3 show that in some cases the obtained metric vanishes only at a sole point and the domain where the metric is positive can be mapped onto the whole 2-plane via an appropriate choice of a pair of conjugate harmonic functions.

Theorem 1 and Lemma 1 imply that any solution to the equation (3.1) provides a certain example of a metric admitting a rational integral. It is an interesting question whether any metric admitting such an integral satisfies the relation (2.16) for a certain solution to (3.1) and a certain pair of conjugate harmonic functions. Another interesting problem is to investigate whether it is possible to provide certain sufficient conditions when the metrics obtained from Theorem 2 are defined and do not vanish on the whole plane. These questions should be investigated in the forthcoming papers.

It is left to notice that the proposed methods also may be useful in construction of metrics admitting rational integrals with a numerator and a denominator of higher degrees.

\vspace{1cm}

Sergei Agapov (corresponding author)

\vspace{2mm}

Novosibirsk State University,

1, Pirogova str., Novosibirsk, 630090, Russia;

\vspace{2mm}

Sobolev Institute of Mathematics SB RAS,

4 Acad. Koptyug avenue, 630090 Novosibirsk Russia.

\vspace{2mm}

agapov.sergey.v@gmail.com, agapov@math.nsc.ru

\vspace{8mm}

Vladislav Shubin

\vspace{2mm}

Novosibirsk State University,

1, Pirogova str., Novosibirsk, 630090, Russia.

\vspace{2mm}

vlad.v.shubin@gmail.com


\vspace{2mm}


\begin{thebibliography}{1}


\bibitem{1} Kozlov V.V.:
{Topological obstacles to the integrability of natural mechanical systems}, Dokl. Akad. Nauk SSSR, {\bf 249}:6 (1979), 1299 -- 1302.

\bibitem{2} Kozlov V.V.:
{Symmetries, topology, and resonances in Hamiltonian mechanics}, Springer-Verlag, Berlin, 1996.

\bibitem{3} Birkhoff G.D.:
{Dynamical systems.} Vol. 9. American Mathematical Society Colloquium Publications, New York (1927).

\bibitem{4} Kolokol'tsov V.N.:
{Geodesic flows on two-dimensional manifolds with an additional first integral that is polynomial in the velocities}, Math. USSR Izv, {\bf 46}:5 (1982), 291 -- 306.

\bibitem{5} Ten V.V.:
{Local integrals of geodesic flows.} Regul. Chaotic Dyn. {\bf 2} (1997), 87 -- 89 [Russian].

\bibitem{6} Bagderina Yu.Yu.:
{Rational integrals of the second degree of two-dimensional geodesic equations}, Sib. Electron. Math. Rep., {\bf 14} (2017), 33 -- 40 [Russian].

\bibitem{7} Collinson C.D.:
{A note on the integrability conditions for the existence of rational first integrals of
the geodesic equations in a Riemannian space.} Gen. Relativity Gravitation, {\bf 18}:2 (1986), 207 -- 214.

\bibitem{8} Collinson C.D., O'Donnell P.J.:
{A class of empty spacetimes admitting a rational first integral of the
geodesic equation.} Gen. Relativity Gravitation, {\bf 24}:4 (1992), 451 -- 455.

\bibitem{9} Combot Th.:
{Rational integrability of trigonometric polynomial potentials on the flat torus.} Regul. Chaotic Dyn., {\bf 22}:4 (2017), 386 -- 397.

\bibitem{10} Heilbronn G.:
{''Integration des equations differentielles ordinaires par la methode de Drach''}, Gauthier-Villars, Paris, 1956.

\bibitem{11} Pavlov M.V., Tsarev S.P.:
{Classical mechanical systems with one-and-a-half degrees of freedom
and Vlasov kinetic equation.} Amer. Math. Soc. Transl., {\bf 234} (2014), 337 -- 371.

\bibitem{12} Maciejewski A.J., Przybylska M.:
{Darboux polynomials and first integrals of natural polynomial Hamiltonian systems.} Phys. Lett. A, {\bf 326}:3-4 (2004), 219 -- 226.

\bibitem{13} Agapov S.V.:
{Rational integrals of a natural mechanical system on the 2-torus.} Sib. Math. Journ. {\bf 61}:2 (2020), 199 -- 207.

\bibitem{14} Agapov S.V.:
{On first integrals of two-dimensional geodesic flows.} Sib. Math. Journ. {\bf 61}:4 (2020), 563 -- 574.

\bibitem{15} Darboux G:
{Lessons on the general theory of surfaces and the geometric applications of infinitesimal calculus.} 1887. Vol. 3.

\bibitem{16} Kozlov V.V.:
{''On Rational Integrals of Geodesic Flows''}, Regul. Chaotic Dyn., {\bf 19}:6 (2014), 601 -- 606.

\bibitem{17} Aoki A., Houri T., Tomoda K.:
{Rational first integrals of geodesic equations and generalised hidden symmetries.} Classical Quantum Gravity, {\bf 33}:19 (2016), 195003, 12 pp.

\bibitem{18} Hietarinta J.:
{New integrable Hamiltonians with transcendental invariants.} Phys. Rev. Lett., {\bf 52}:1057 (1984).

\bibitem{19} Perelomov A.M.:
{Integrable Systems of Classical Mechanics and Lie Algebras.} Birkhauser Verlag Basel, 1990.

\bibitem{20} Abramowitz M., Stegun I. A.:
{Handbook of Mathematical Functions with Formulas, Graphs, and Mathematical Tables.} National Bureau of Standards Applied Mathematics Series 55. Tenth Printing. (1972).

\bibitem{21} Bitsadze A.V.:
{Osnovy teorii analiticheskih funkciy compleksnogo peremennogo.} Nauka, Moscow, 1969 [Russian].

\end{thebibliography}
\end{document}